\documentclass[12pt,a4paper]{amsart}
\usepackage{amssymb,amscd}
\title[Natural differential operations]{Natural differential operations on manifolds: an algebraic approach}
\thanks{This research was supported by the RFBR grant 05--01--00988.}
\date{July 4, 2006}
\author{P. I. Katsylo}
\address{Independent University of Moscow\\
Bolshoi Vlasievskii 11, 119002, Moscow, Russia}
\email{p57taras@yandex.ru}
\author{D. A. Timashev}
\address{Moscow State University, Faculty of Mechanics and Mathematics\\
Leninskie gory, 119992, Moscow, Russia}
\email{timashev@mech.math.msu.su}
\urladdr{http://mech.math.msu.su/department/algebra/staff/timashev}

\newcommand{\KK}{\mathbb{K}}
\newcommand{\RR}{\mathbb{R}}
\newcommand{\CC}{\mathbb{C}}
\newcommand{\ZZ}{\mathbb{Z}}
\newcommand{\M}{M}
\newcommand{\n}{n}
\newcommand{\kk}{k}
\newcommand{\pt}{z}
\newcommand{\x}[2][]{x\ifx#1\empty\else_{#1}\fi\ifx#2\empty\else^{#2}\fi}
\newcommand{\e}[1]{e\ifx#1\empty\else_{#1}\fi}
\newcommand{\y}[1]{y\ifx#1\empty\else_{#1}\fi}
\newcommand{\g}[1][]{g\ifx#1\empty\else_{#1}\fi}
\newcommand{\dg}[1][]{\xi\ifx#1\empty\else_{#1}\fi}
\newcommand{\f}[1][]{f\ifx#1\empty\else_{#1}\fi}
\newcommand{\vv}[1][]{v\ifx#1\empty\else_{#1}\fi}
\newcommand{\J}[3][]{J\ifx#1\empty\else_{#1}\fi^{#2}#3}
\newcommand{\JP}[2][\kk]{J_{#2}^{(#1)}}
\newcommand{\GL}[1]{\mathrm{GL}_{#1}}
\newcommand{\gl}[1]{\mathfrak{gl}_{#1}}
\newcommand{\JGL}[2][\kk]{\mathrm{GL}_{#2}^{(#1)}}
\newcommand{\jgl}[2][\kk]{\mathfrak{gl}_{#2}^{(#1)}}
\newcommand{\NGL}[2][\kk]{\mathrm{NGL}_{#2}^{(#1)}}
\newcommand{\ngl}[2][\kk]{\mathfrak{ngl}_{#2}^{(#1)}}
\newcommand{\SL}[1]{\mathrm{SL}_{#1}}
\newcommand{\Or}[1]{\mathrm{O}_{#1}}
\newcommand{\SO}[1]{\mathrm{SO}_{#1}}
\newcommand{\Sp}[1]{\mathrm{Sp}_{#1}}
\newcommand{\spl}[1]{\mathfrak{sp}_{#1}}
\newcommand{\JSp}[2][\kk]{\mathrm{Sp}_{#2}^{(#1)}}
\newcommand{\NSp}[2][\kk]{\mathrm{NSp}_{#2}^{(#1)}}
\newcommand{\jsp}[2][\kk]{\mathfrak{sp}_{#2}^{(#1)}}
\newcommand{\nsp}[2][\kk]{\mathfrak{nsp}_{#2}^{(#1)}}
\newcommand{\dual}[1]{#1^{*}}
\newcommand{\metric}{\omega}
\newcommand{\Sym}[1][\bullet]{\mathrm{S}^{#1}}
\newcommand{\Ext}[1][\bullet]{{\textstyle\bigwedge^{#1}}}
\newcommand{\Sch}[1]{\mathbb{S}^{#1}}
\newcommand{\Symm}[1][]{\mathop{\mathrm{Sym}}\nolimits_{#1}}
\newcommand{\Alt}[1][]{\mathop{\mathrm{Alt}}\nolimits_{#1}}
\newcommand{\Ker}{\mathop{\mathrm{Ker}}}
\newcommand{\rk}{\mathop{\mathrm{rk}}}
\newcommand{\R}{R}
\newcommand{\dR}{\rho}
\newcommand{\act}{\gamma}
\newcommand{\contr}{\sigma}
\newcommand{\Rep}[2][\kk]{\mathop{\mathrm{Fr}}\nolimits^{#1}(#2)}
\newcommand{\SRep}[2][\kk]{\mathop{\mathrm{SFr}}\nolimits^{#1}(#2)}
\newcommand{\Ass}[3][\kk]{\Rep[#1]{#2}\gtimes{\JGL[#1]\n}#3}
\newcommand{\gtimes}[1]{\mathbin{\times^{#1}}}
\newcommand{\C}{\mathcal{C}}
\newcommand{\F}{\mathcal{F}}
\newcommand{\V}{\mathcal{V}}
\newcommand{\W}{\mathcal{W}}
\newcommand{\OO}{\mathcal{O}}
\newcommand{\T}[2][]{\mathcal{T}\ifx#2\empty\else^{\ifx#1\empty\else#1,\fi#2}\fi}
\newcommand{\DF}[1][\bullet]{\Omega^{#1}}
\newcommand{\JF}[2][\kk]{{#2^{(#1)}}}
\newcommand{\Ff}{\mathfrak{F}}
\newcommand{\Man}[1]{\mathcal{M}an_{#1}}
\newcommand{\Fib}[1]{\mathcal{F}ib_{#1}}
\newcommand{\DD}{D}
\newcommand{\dd}{\delta}
\newcommand{\Mor}[3][]{\mathop{\mathrm{Mor}}\nolimits_{#1}(#2,#3)}
\newcommand{\D}[3][]{\mathcal{D}\ifx#1\empty\else^{#1}\fi(#2,#3)}
\newcommand{\Sec}[2][]{\Gamma(\ifx#1\empty\else#1,\fi#2)}

\newcommand{\reg}{\text{reg}}
\newcommand{\GS}{\mathcal{A}}
\newcommand{\Gs}{A}
\newcommand{\gs}{\alpha}
\renewcommand{\[}{\mathopen{[\mspace{-1mu}[}}
\renewcommand{\]}{\mathclose{]\mspace{-1mu}]}}
\newcommand{\Inv}{\mathcal{M}}
\newcommand{\h}{\varepsilon}

\renewcommand{\le}{\leqslant}
\renewcommand{\ge}{\geqslant}

\newtheorem{theorem}{Theorem}
\newtheorem{lemma}{Lemma}
\newtheorem*{corollary}{Corollary}
\theoremstyle{definition}
\newtheorem{definition}{Definition}
\newtheorem{example}{Example}
\theoremstyle{remark}
\newtheorem*{remark}{Remark}
\newtheoremstyle{claim}{}{}{\itshape}{0pt}{\bfseries}{.}{.5em}{\thmnote{#3}}
\theoremstyle{claim}
\newtheorem*{claim}{}
\newtheoremstyle{concept}{}{}{}{0pt}{\bfseries}{.}{.5em}{\thmnote{#3}}
\theoremstyle{concept}
\newtheorem*{concept}{}

\begin{document}

\begin{abstract}
We consider natural algebraic differential operations acting on
geometric quantities over smooth manifolds. We introduce a method
of study and classification of such operations, called
IT-reduction. It reduces the study of natural operations to the
study of polynomial maps between (vector) spaces of jets which are
equivariant with respect to certain algebraic groups. Using the
IT-reduction, we obtain short and conceptual proofs of some known
results on the classification of certain natural operations (the
Schouten theorem, etc) together with new results including the
non-existence of a universal deformation quantization on Poisson
manifolds.
\end{abstract}

\maketitle

\section*{Introduction}

In differential geometry, there are many nontrivial formul{\ae}
and theorems based on local calculations. Typical examples are:
the Bianchi identities, the Gilkey theorem~\cite{Gilkey}, the
Weitzenb{\"o}ck formula \cite[Ch.\,1, I]{Einstein}. Such theorems
and formul{\ae} of purely local nature lie in the basis of
differential geometry. Moreover, it often happens that a discovery
of some local formula allows to solve an important problem. For
example, the Gilkey theorem led to a new proof of the index
theorem \cite{index}, \cite{Gilkey&index}. Thus it is tempting to
look for a universal approach to local problems of differential
geometry.

One of possible approaches, called \emph{formal geometry}, was
suggested by I.~M.~Gelfand and D.~A.~Kazhdan
in~1971~\cite{formal}. Approximately at the same time,
E.~B.~Vinberg observed that using simple arguments from the
representation theory of algebraic groups and invariant theory
makes local calculations in differential geometry much easier and
more intelligible. (This was one of initial motivations for the
study of invariant theory at the Moscow school of invariant theory
leaded by Vinberg.) This approach reduces local problems of
differential geometry to problems in the invariant theory of
\emph{finite-dimensional representations of linear algebraic
groups}. We call it the invariant-theoretic reduction
(\emph{IT-reduction}).

A simplified version of the IT-reduction method is exposed
in~\cite{curv}. In this paper we introduce this method in maximal
generality which is necessary for applications. The idea of local
study of geometric quantities and natural differential operations
between them by considering jet spaces and the actions of
coordinate transformations on them was exploited by many
researchers, see e.g.\ the monograph~\cite{nat.op}. We develop an
instrumental approach concentrating on effective
rep\-re\-sen\-ta\-tion- and invariant-theoretic methods for
solving concrete local problems of differential geometry. For
instance, in \cite{nat.op} the concept of natural operations is
developed in very general context and much effort is put on
proving that certain natural operations are of finite order, using
Peetre-type theorems, while we impose the finite order assumption
from the very beginning and consider mainly differential
operations given by algebraic formul{\ae} believing that this case
is most interesting and essential in applications. We illustrate
the IT-reduction method by solving several local problems. Some of
them were previously solved by other authors using ad hoc methods,
and other results are new.

Now we briefly describe the content of the paper.

In~\S\ref{basic} we introduce basic notions of geometric objects
and quantities, natural bundles and differential operations.
Essentially, geometric quantities (e.g.~tensor fields) are
sections of fibre bundles associated with coframe bundles and
natural differential operations act on geometric quantities by
universal differential formul{\ae} that are invariant under
coordinate transformations. Morally, such operations, like
exterior differential or the curvature of a Riemannian metric,
should have intrinsic ``physical meaning'' since they do not
depend on a chosen frame of reference. In this section, we explain
the IT-reduction method. Also, we recall some basic facts from the
representation theory of classical groups, which are used in
computations.

Then we concentrate on polynomial natural differential operations
acting on tensor fields. Starting with some simple reductions and
general finiteness results in~\S\ref{finite}, we derive in
\S\ref{linear} the classification of natural linear differential
operations, which goes back to Schouten. In \S\ref{symplectic} we
use the IT-reduction to extend this result to manifolds with
additional symplectic structure, cf.~\cite{irrep(ham)}. This
requires an extension of the notion of a natural differential
operation to manifolds equipped with an additional structure,
see~\ref{geom.str}. Finally, we prove in~\S\ref{def.quant} that
there exists no universal formula for deformation quantization on
Poisson manifolds which is invariant under coordinate
transformations.

\subsubsection*{Convention}

In formul{\ae} of tensor calculus, we systematically use the
Einstein summation rule, i.e., assume by default the summation in
each pair of coinciding upper and lower indices which run from 1
up to the dimension of a manifold.

\section{Preliminaries}
\label{basic}

We work over smooth real or complex manifolds. However our
considerations will be purely algebraic and the ground field $\KK$
will not play any essential r\^{o}le. So we consider smooth
$\KK$-manifolds making no distinction between the cases $\KK=\RR$
(differential geometry) and $\KK=\CC$ (complex analytic geometry).

\subsection{Geometric quantities}
\label{geom.quant}

The concept of a geometric quantity goes back to Riemann, see
\cite[Ch.\,6, \S1]{diff.geom}. Loosely speaking, a geometric
quantity is a scalar value or a tuple of scalars associated with
each point of a manifold in a way depending on chosen local
coordinates which transforms under a coordinate change in a
regular way depending only on initial values and partial
derivatives, up to a certain order, of new coordinates with
respect to old ones. Natural examples are tensor fields (order~1)
and connections (order~2). A rigorous definition of a geometric
quantity can be given in several ways.

The classical analytic definition sounds as follows. Let $\M$ be a
manifold of dimension~$\n$. A \emph{geometric quantity} is a
function $\f[\alpha]=\f[\alpha](\x[\alpha]{})$ of local
coordinates $\x[\alpha]{}=(\x[\alpha]1,\dots,\x[\alpha]\n)$ on
$\M$ which takes values in a space (vector space or, more
generally, manifold) $F$ and transforms under a coordinate change
$\x[\alpha]{}\to\x[\beta]{}$ in the following way:
\begin{equation}\label{anal.def}
\f[\beta]= \Phi\left(\left\{
\frac{\partial^l\x[\beta]i}{\partial(\x[\alpha]1)^{l_1}\cdots
\partial(\x[\alpha]\n)^{l_{\n}}}\right\}_
{\substack{1\le l\le\kk,\ 1\le i\le\n,\\l_1+\dots+l_{\n}=l}},\
\f[\alpha]\right),
\end{equation}
where $\Phi$ is a differentiable map. In other words,
$\f[\beta]$~depends only on $\f[\alpha]$ and on the $\kk$-jet
$\J\kk{\g[\beta\alpha]}$ of the coordinate transformation
$\x[\beta]{}=\g[\beta\alpha](\x[\alpha]{})$.

In order to reformulate this definition in modern terms, consider
the group $\JGL\n$ of $\kk$-jets of local diffeomorphisms
$\KK^{\n}\to\KK^{\n}$ at~$0$. Elements of $\JGL\n$ are represented
in the form:
\begin{multline*}
\x{}\mapsto\g(\x{})=\g[1](\x{})+\g[2](\x{},\x{})+\dots+\g[\kk](\x{},\dots,\x{}),\\
\x{}=(\x1,\dots,\x\n),\quad
\g[l]\in\Sym[l](\KK^{\n})^*\otimes\KK^{\n}\
(l=1,\dots,\kk),\quad\det\g[1]\ne0.
\end{multline*}
$\JGL\n$~is a linear algebraic group isomorphic to the
automorphism group of the truncated polynomial algebra
\begin{equation*}
\JP\n=\KK[\x1,\dots,\x\n]/(\x1,\dots,\x\n)^{\kk+1}.
\end{equation*}
Its unipotent radical $\NGL\n$ is defined by the equation
$\g[1](\x{})=\x{}$ and the subgroup $\GL{\n}$ of linear
transformations is a Levi subgroup.

The Lie algebra $\jgl\n$ of $\JGL\n$ is identified with the space
of polynomial vector fields of degree $\le\kk$ vanishing at~0,
which are represented by polynomial maps $\x{}\mapsto\g(\x{})$ as
above, but without the restriction $\det\g[1]\ne0$. The Lie
algebra $\ngl\n$ of $\NGL\n$ is distinguished by $\g[1]=0$.

It is easy to see that the map $\Phi:\JGL\n\times F\to F$
in~(\ref{anal.def}) defines an action of $\JGL\n$ on~$F$. In most
applications, $F$~is a vector or affine space, or an open subset
in such a space, and the action $\JGL\n:F$ is a rational linear or
affine representation.

This consideration leads to a geometric reformulation of the above
analytic definition of geometric quantities.

Let $\Rep\M$ denote the coframe bundle of order $\kk$ on~$\M$
\cite[Ch.\,6, 1.2]{diff.geom}. The fibre of $\Rep\M\to\M$ over
$\pt\in\M$ consists of \emph{coframes of order $\kk$} at~$\pt$,
i.e., $\kk$-jets of coordinate systems $\x{}=(\x1,\dots,\x\n)$ in
a neighborhood of $\pt$ with $\x{}(\pt)=0$. $\Rep\M$~is a
principal bundle with respect to the natural action of~$\JGL\n$.
Instead of coframes, one may consider \emph{frames of
order~$\kk$}, which are $\kk$-jets of inverse coordinate maps
sending a neighborhood of~$0\in\KK^{\n}$ onto a neighborhood
of~$\pt\in\M$. Frames of order~1 are determined by fixing a basis
of~$T_{\pt}\M$ (~=~the image of the standard basis of
$T_0\KK^{\n}=\KK^{\n}$), i.e., a usual frame on~$\M$. The bundles
of frames and coframes of order $\kk$ are canonically isomorphic.

\begin{concept}[Geometric definition]
Suppose that $F$ is a manifold equipped with a differentiable
action of~$\JGL\n$. The associated fibre bundle
\begin{equation*}
\F=\Ass{\M}F:=(\Rep\M\times F)/\JGL\n
\end{equation*}
is said to be the \emph{space of geometric objects of type $F$}
on~$\M$. If the action $\JGL\n:F$ is not reduced to the action of
the quotient group~$\JGL[\kk-1]\n$, then we say that geometric
objects of type $F$ have \emph{order~$\kk$}.

Spaces of geometric objects are also called \emph{natural
bundles}.

A \emph{geometric quantity of type $F$} is a section of the
natural bundle $\F\to\M=\Rep\M/\JGL\n$. The set of geometric
quantities is denoted by $\Sec\F=\Sec[\M]\F$.
\end{concept}

\begin{example}
Let $F=(\KK^{\n})^{\otimes p}\otimes(\KK^{\n}{}^{*})^{\otimes q}$,
with the natural linear action of $\GL{\n}=\JGL[1]\n$. Then
$\F=\T[p]q$ is the tensor bundle of type $(p,q)$.

Generally, every natural vector bundle $\F$ of order~1, with an
additional requirement that the representation $\GL\n:F$ be
rational, embeds as a subbundle into a product of $\T[p]q$'s. For
this reason, we call such $\F$ \emph{tensor bundles}.
\end{example}

\begin{example}[{\cite[Ch.\,6, \S4]{diff.geom}}]
The $l$-jets of local sections of $\F\to\M$ form the \emph{$l$-jet
bundle} $\JF[l]\F$ of order~$k+l$. Every section $\f:\M\to\F$
defines a section $\J{l}\f:\M\to\JF[l]\F$ whose value at
$\pt\in\M$ is the $l$-jet $\J[\pt]l\f$ of $\f$ at~$\pt$. The
action of $\JGL[\kk+l]\n$ on the typical fibre $\JF[l]F$ of
$\JF[l]\F$ is derived from the transition rule for local
coordinate presentations of a section of~$\F$. It is given by the
formula
\begin{equation}\label{jet.act}
\J[0]{\kk+l}\g\cdot\J[0]l\f=\J[0]l\left(\J[\g^{-1}(\x{})]\kk\g\cdot
\f\bigl(\g^{-1}(\x{})\bigr)\right),
\end{equation}
for every local diffeomorphism $\g$ of $\KK^{\n}$ with $\g(0)=0$
and local section $\f$ of $\KK^{\n}\times F\to \KK^{\n}$ defined
in a neighborhood of~0.

If $F$ is a vector space, then $\JF[l]F=\JP[l]\n\otimes F$.
(\emph{Caution}: the $\JGL[\kk+l]\n$-ac\-tion on $\JF[l]F$ is
\emph{not} given by the tensor product of linear representations.)
\end{example}


For the sequel, we need a more explicit description of how jets of
diffeomorphisms act on jets of tensor fields.

\begin{lemma}\label{tensor.jet}
Let $\R:\GL\n\to\GL{}(F)$ be a rational representation and
$\dR:\gl\n\to\gl{}(F)$ the corresponding Lie algebra
representation. The natural actions $\JGL[\kk+1]\n:\JF{F}$ and
$\jgl[\kk+1]\n:\JF{F}$ are given by the formul{\ae}:
\begin{multline*}
\begin{aligned}
(\g\cdot\f)(\x{})&=\R\left(\sum_{l=1}^{\kk+1}l\cdot
\g[l]\bigl(\g^{-1}(\x{}),\dots,\g^{-1}(\x{}),\cdot\bigr)\right)
\f\bigl(\g^{-1}(\x{})\bigr),\\
(\dg\cdot\f)(\x{})&=\sum_{l=1}^{\kk+1}l\cdot
\dR\bigl(\dg[l](\x{},\dots,\x{},\cdot)\bigr)\f(\x{})-
\sum_{l=1}^{\kk}l\cdot\f[l]\bigl(\x{},\dots,\x{},\dg(\x{})\bigr),
\end{aligned}\\
\forall\g\in\JGL[\kk+1]\n,\ \dg\in\jgl[\kk+1]\n,\ \f\in F.
\end{multline*}
Here $\g[l],\dg[l]\in\Sym[l]\KK^{\n*}\otimes\KK^{\n}$ and
$\f[l]\in\Sym[l]\KK^{\n*}\otimes F$ are the homogeneous components
of $\g,\dg,\f$, and the r.h.s.\ are interpreted as follows. The
arguments of $\R,\dR$ are polynomial functions in $\x{}$ taking
values in $\GL\n$, resp.~in~$\gl\n$, at $\x{}$ from a neighborhood
of~0. The r.h.s.\ are expanded in Taylor series in~$\x{}$ and
truncated from the order~$\kk+1$.
\end{lemma}

\begin{corollary}
\begin{equation}\label{Lie.act}
(\dg\cdot\f)_{\kk}=\sum_{l=1}^{\kk+1}l\act(\dg[l]\otimes
\f[\kk+1-l])-(\kk+1-l)\contr(\dg[l]\otimes\f[\kk+1-l]),
\end{equation}
where $\act,\contr:\Sym[l]\KK^{\n*}\otimes\KK^{\n}\otimes
\Sym[\kk+1-l]\KK^{\n*}\otimes F\to\Sym[\kk]\KK^{\n*}\otimes F$ are
linear maps defined as follows. To compute~$\act$, consider
$\Sym[l]\KK^{\n*}\otimes\KK^{\n}$ as a subspace of
$\Sym[l-1]\KK^{\n*}\otimes\gl\n$, apply $\gl\n$ to $F$ via~$\dR$,
and conclude by the multiplication
$\Sym[l-1]\KK^{\n*}\otimes\Sym[\kk+1-l]\KK^{\n*}\to\Sym[\kk]\KK^{\n*}$.
The map $\contr$ is the contraction of $\KK^{\n}$ with
$\Sym[\kk+1-l]\KK^{\n*}$ followed by the multiplication
$\Sym[l]\KK^{\n*}\otimes\Sym[\kk-l]\KK^{\n*}\to\Sym[\kk]\KK^{\n*}$.
\end{corollary}

The proof is a straightforward calculation based
on~(\ref{jet.act}).

Note that the space of geometric objects of given type is
naturally defined over any $\n$-dimensional manifold. This
observation leads to a reformulation of the definition of
geometric objects in the categorical language~\cite{geom.obj}.

Let $\Man\n$ denote the category of $\n$-dimensional manifolds
where the morphisms are open embeddings. Let $\Fib\n$ denote the
category of fibre bundles over $\n$-dimensional manifolds, the
morphisms being differentiable maps of bundles covering the
morphisms of their bases in~$\Man\n$.

\begin{concept}[Categorical definition]
A \emph{type of geometric objects} is a functor
$\Ff:\Man\n\rightsquigarrow\Fib\n$ such that $\Ff(\M)$ is a bundle
over $\M$ for any $\n$-di\-men\-sional manifold $\M$ and
$\Ff(\M')$ is the restriction of $\Ff(\M)$ for any open
submanifold $\M'\subset\M$, the morphism $\Ff(\M')\to\Ff(\M)$
induced by $\M'\hookrightarrow\M$ being the inclusion.
\end{concept}

It is easy to see that all $\F=\Ff(\M)$ have one and the same
typical fibre~$F$. Every local diffeomorphism
$\g[\beta\alpha]:U_{\alpha}\to U_{\beta}$ between two sufficiently
small neighborhoods of $0\in\KK^n$ with $\g[\beta\alpha](0)=0$
induces a diffeomorphism
$\Ff(\g[\beta\alpha]):\Ff(U_{\alpha})\simeq U_{\alpha}\times
F\to\Ff(U_{\beta})\simeq U_{\beta}\times F$. These
$\Ff(\g[\beta\alpha])$ are compatible with shrinking of
$U_{\alpha},U_{\beta}$ and thus induce transformations of $F$
(~=~the fibre at~$0$) depending only on the germs of
$\g[\beta\alpha]$ at~$0$. Palais and Chuu-Lian Terng
\cite{geom.obj} proved that in fact these transformations of $F$
depend only on $\J\kk\g[\beta\alpha]$ for sufficiently
large~$\kk$. This yields an action $\JGL\n:F$ and isomorphism
$\Ff(\M)\simeq\Ass{\M}F$.

\subsection{Differential operations}

Differential operations act on geometric quantities. Given two
natural bundles $\V,\W\to\M$, a differential operation $\DD$ from
$\V$ to $\W$ transforms (local) sections of $\V$ into those of
$\W$ according to a formula of the following kind, in local
coordinates:
\begin{equation}\label{diff.op}
(\DD\vv)^p=\dd^p\left(\left\{\x{i},\vv^q,\partial_1^{l_1}\!\cdots
\partial_{\n}^{l_{\n}}\vv^q\right\}_{\substack{1\le i\le\n,\ 1\le q\le\dim V,\\
1\le l_1+\dots+l_{\n}\le\kk}}\right),\quad\forall\vv\in\Sec\V,
\end{equation}
where $\dd^p$ are differentiable functions ($p=1,\dots,\dim W$),
$V,W$~are typical fibres of~$\V,\W$, and
$\partial_j=\partial/\partial\x{j}$ ($j=1,\dots,\n$).

In other words, a \emph{differential operation of order $\le\kk$}
is a map $\DD:\Sec\V\to\Sec\W$ induced by a morphism $\JF\V\to\W$
(denoted by the same letter) so that
$\DD\vv(\pt)=\DD\bigl(\J[\pt]\kk\vv\bigr)$, $\forall\vv\in\Sec\V$,
$\pt\in\M$ \cite[Ch.\,6, 4.6]{diff.geom}. (The order is exactly
$\kk$ if $\DD$ does not factor through the canonical projection
$\JF\V\to\JF[\kk-1]\V$.)

\begin{remark}
$\DD$~produces a series of differential operations
$\JF[\kk+l]\V\to\JF[l]\W$ (denoted by the same letter) in an
obvious way.
\end{remark}

\begin{definition}\label{alg.op}
Suppose that $V,W$ are open invariant subsets in vector or affine
spaces equipped with rational representations of some~$\JGL[l]\n$.
A differential operation $\DD:\JF\V\to\W$ is \emph{algebraic} if
the maps of fibres $\dd_{\x{}}:\JF{V}\to W$ are algebraic
morphisms of bounded degree. In other words,
$\dd^p$~in~(\ref{diff.op}) are rational functions in
$\vv^q,\partial_1^{l_1}\!\cdots\partial_{\n}^{l_{\n}}\vv^q$ whose
coefficients are differentiable functions in~$x$, denominators
depend only on~$\vv^q$, and the degrees of numerators are bounded
on~$\M$.
\end{definition}

\begin{remark}
One may extend Definition~\ref{alg.op} replacing typical fibres
$V,W$ by more general algebraic varieties with algebraic
$\JGL[l]\n$-actions. However our formulation is sufficient for
many applications.
\end{remark}

We focus our attention at algebraic differential operations as the
most customary case.

Algebraic differential operations may be regarded as geometric
objects, too. In the simplest case where $V,W$ are vector spaces,
all $\dd_{\x{}}$~are polynomial maps of degree~$\le d$. The set
$\Mor[d]{\JF{V}}W$ of polynomial maps $\JF{V}\to W$ of degree $\le
d$ is a vector space with the natural action of $\JGL[\kk+l]\n$ by
conjugation. It is easy to see that $\DD$ is nothing but a section
of $\Ass[\kk+l]\M{\Mor[d]{\JF{V}}W}$. The general case is handled
in the same way if one restricts from above the degrees of
denominators, too.

\begin{example}
Let $\V,\W$ be vector bundles. Linear differential operators of
order $\le\kk$ from $\V$ to $\W$ are geometric quantities taking
values in $\D[\kk]\V\W=\Ass[\kk+l]\M{\bigl(\JF{V}^{*}\otimes
W\bigr)}$.
\end{example}

\begin{example}
Let $\T{}=\T[1]0$ be the tangent bundle of~$\M$. Consider the
subbundle $\C\subset\D[1]{\T{}}{\T{}\otimes\T*}$ consisting of
homomorphisms $\nabla$ splitting the natural exact sequence:
\begin{equation*}
0\longrightarrow\T{}\otimes\T*
\underset{\textstyle\underset\nabla\dashleftarrow}\longrightarrow
\JF[1]{\T{}}\longrightarrow\T{}\longrightarrow0.
\end{equation*}
Sections $\nabla\in\Sec\C$ act on vector fields as covariant
derivations: in local coordinates, given a vector field
$\xi(\x{})=\xi_0+\xi_j\x{j}+\cdots$ in a neighborhood of
$\pt\in\M$, we have
\begin{equation*}
\nabla\xi(\pt)=\nabla(\xi_0^i\partial_i)+\nabla(\xi_j^i\x{j}\partial_i)=
\Gamma^k_{ij}\xi_0^i\;\partial_k\otimes d\x{j}+
\xi_j^i\;\partial_i\otimes d\x{j}.
\end{equation*}
Hence geometric quantities with values in $\C$ are linear
connections on~$\M$. They are affine geometric objects of order~2.
\end{example}

\subsection{Natural operations}
\label{nat.op}

Natural differential operations on geometric quantities are
distinguished by the property that their coordinate expression is
one and the same for any choice of local coordinates. This
property may be reformulated as follows.

\begin{definition}
A differential operation $\DD:\JF\V\to\W$ is called \emph{natural}
if the respective map $\dd:\JF{V}\to W$ of fibres at $\pt\in\M$
does not depend on~$\pt$ and is $\JGL[\kk+l]\n$-equivariant
(assuming that $\V,\W$ have order $\le l$).
\end{definition}

Since a natural differential operation $\DD$ is uniquely
determined by the map of typical fibres~$\dd$, it follows that
$\DD$ is naturally defined on geometric quantities of given type
over any $\n$-dimensional manifold. This observation leads to a
functorial point of view on natural operations.

For each type of geometric objects~$\Ff$, consider the respective
functor of geometric quantities
$\M\rightsquigarrow\Sec[\M]\F=\Sec\F$, which associates with an
$\n$-manifold $\M$ the space of sections of $\F=\Ff(\M)\to\M$. It
becomes a contravariant functor from $\Man\n$ to topological
spaces, if we equip $\Sec\F$ with the topology of locally uniform
convergence of sections and all their partial derivatives.

A natural operation $\DD:\JF\V\to\W$ induces a natural
transformation of functors: there is a commutative square
\begin{equation*}
\begin{CD}
\Sec[\M]\V  & @>\DD>> & \Sec[\M]\W  \\
  @VVV      &         &    @VVV     \\
\Sec[\M']\V & @>\DD>> & \Sec[\M']\W,
\end{CD}
\end{equation*}
for every open embedding $\M'\hookrightarrow\M$. Conversely, under
certain conditions a natural transformation of geometric
quantities is given by a natural differential operation, by
Peetre-type theorems \cite[Ch.\,5]{nat.op}.

Algebraic natural operations are given by everywhere defined
rational maps $\dd:\JF{V}\to W$ that are equivariant with respect
to the action of an appropriate algebraic group~$\JGL[\kk+l]\n$.
Thus the study of such operations is a \emph{purely algebraic
problem} belonging to the representation theory of algebraic
groups and invariant theory. For this reason, we call our approach
to algebraic differential operations the \emph{IT-reduction}.
(IT~stands for ``invariant-theoretic''.)

In this paper, we consider natural algebraic differential
operations on tensor bundles. By the above discussion, they are in
a bijective correspondence with polynomial (or, more generally,
rational) maps $\dd:\JP\n\otimes V\to W$ that are
$\GL\n$-equivariant and $\NGL[\kk+1]\n$-invariant, where $V,W$ are
rational representations of~$\GL\n$ and $\JGL[\kk+1]\n$ acts on
$\JP\n\otimes V=\JF{V}$ in the natural way
(cf.~Lemma~\ref{tensor.jet}).

\begin{example}\label{ext.diff}
A classical example of a natural (algebraic) differential
operation is the exterior differential
$d:\Sec{\DF[m]}\to\Sec{\DF[m+1]}$, where $\DF[m]=\Ext[m]\T*$ is
the bundle of exterior $m$-forms. The respective map of fibres
\begin{equation*}
\dd:\JP[1]\n\otimes\Ext[m]\KK^{\n*}\to\KK^{\n*}\otimes\Ext[m]\KK^{\n*}\to\Ext[m+1]\KK^{\n*}
\end{equation*}
is the canonical projection (ignoring the constant term) followed
by the alternation:
\begin{multline*}
\dd(\y0\otimes\y1\wedge\dots\wedge\y{m})=
d(\y0\:d\y1\wedge\dots\wedge d\y{m})=
\y0\wedge\y1\wedge\dots\wedge\y{m},\\
\forall\y0,\dots,\y{m}\in\KK^{\n*}.
\end{multline*}
\end{example}

\begin{example}
Lie derivative on a tensor bundle $\V$ may be regarded as a
natural bilinear operation of order~1 from $\T{}\times\V$ to~$\V$.
The respective map of fibres is
\begin{align*}
\dd:\bigl(\KK^{\n}\oplus(\KK^{\n*}\otimes\KK^{\n})\bigr)\times
\bigl(V\oplus(\KK^{\n*}\otimes V)\bigr)&\to V, \\
(\xi_0+\xi_1,\vv[0]+\vv[1])&\mapsto
\dR(\xi_1)\vv[0]-\contr(\xi_0\otimes\vv[1]),
\end{align*}
where $\dR$ is the tensor representation of
$\gl\n\simeq\KK^{\n*}\otimes\KK^{\n}$ in the typical fibre~$V$,
and $\sigma$ denotes the contraction of $\KK^{\n}$
with~$\KK^{\n*}$.
\end{example}

Other examples are curvatures of Riemannian metrics, etc.

\subsection{Classical groups: representations and invariants}

We recall some basic facts about rational representations of
classical linear groups and classical invariant theory, which we
use in the sequel. Our basic references will be \cite{rep.th},
\cite{inv.th}.

A \emph{partition} is a weakly decreasing sequence of non-negative
integers $\lambda=(\lambda_1,\dots,\lambda_{\n})$,
$\lambda_1\ge\dots\ge\lambda_{\n}\ge0$. The \emph{length} of
$\lambda$ is the number of $\lambda_i\ne0$, and
$|\lambda|=\lambda_1+\dots+\lambda_{\n}$ is the number partitioned
by~$\lambda$. Fragments of the form $d,\dots,d$ ($s$~times) in
$\lambda$ are often written as~$d^s$. Omitting the non-negativity
condition $\lambda_{\n}\ge0$ yields the definition of a
\emph{virtual partition}.

Let $\Sch\lambda$ denote the Schur functor corresponding to a
partition~$\lambda$. To any vector space~$V$, it relates a
subspace $\Sch\lambda{V}\subset V^{\otimes|\lambda|}$ constructed
as follows. We may assume that the tensor factors are indexed by
the boxes of the Young diagram corresponding to~$\lambda$. Then
$\Sch\lambda{V}$ is obtained from $V^{\otimes|\lambda|}$ by
applying first the symmetrization in each row of the Young
diagram, denoted~$\Symm[\lambda]$, and then the alternation in
each column~$\Alt[\lambda]$.

$\Sch\lambda{V}$~is an irreducible polynomial $\GL{}(V)$-module
spanned by
\begin{equation*}
\Alt[\lambda]\left(\vv[1]^{\otimes\lambda_1}\otimes\dots\otimes
\vv[\n]^{\otimes\lambda_{\n}}\right),\qquad
\vv[1],\dots,\vv[\n]\in V.
\end{equation*}
Every rational $\GL{}(V)$-module decomposes into a direct sum of
irreducible submodules isomorphic to
$\Sch\lambda{V}\otimes\det^d$, $d\in\ZZ$.

In our considerations it will be convenient to realize the
irreducible rational $\GL\n$-modules as
$\Sch\lambda\KK^{\n*}\otimes\det^d$. Such a module is determined,
up to isomorphism, by a virtual partition
$\overline\lambda=(\overline\lambda_1,\dots,\overline\lambda_{\n})$,
$\overline\lambda_i=\lambda_i-d$. It contains unique, up to
proportionality, eigenvectors with respect to the mutually
opposite lower- and upper-triangular Borel subgroups
$B^-,B^+\subset\GL\n$, namely
\begin{align*}
\vv[\lambda]^-&=
\Alt[\lambda]\left((\x1)^{\otimes\lambda_1}\otimes\dots\otimes(\x\n)^{\otimes\lambda_{\n}}\right),\\
\vv[\lambda]^+&=
\Alt[\lambda]\left((\x\n)^{\otimes\lambda_1}\otimes\dots\otimes(\x1)^{\otimes\lambda_{\n}}\right),
\end{align*}
called \emph{lowest}, resp.~\emph{highest}, \emph{weight vectors}.
Note that $\vv[\lambda]^{\pm}$ generates
$\Sch\lambda\KK^{\n*}\otimes\det^d$ as a $B^{\mp}$-module.

Similarly, the irreducible representations of $\SL\n$ are realized
in~$\Sch\lambda\KK^{\n*}$, $\lambda_{\n}=0$. The irreducible
representations of $\Sp\n$ ($\n$~even) are parameterized by
partitions $\lambda$ of length $l\le\n/2$ and realized in the
subspaces
$\Sch{\langle\lambda\rangle}\KK^{\n*}\subset\Sch\lambda\KK^{\n*}$
spanned by
\begin{equation*}
\Alt[\lambda]\bigl(\y1^{\otimes\lambda_1}\otimes\dots\otimes
\y{l}^{\otimes\lambda_l}\bigr),
\end{equation*}
where $\y1,\dots,\y{l}$ span an isotropic subspace in~$\KK^{\n*}$.
$B^{\pm}\cap\Sp\n$~are mutually opposite Borel subgroups
in~$\Sp\n$, with highest/lowest vectors
$\vv[\lambda]^{\pm}\in\Sch{\langle\lambda\rangle}\KK^{\n*}$,
provided that the symplectic form has secondary-diagonal matrix.

There are effective formul{\ae} for decomposing certain tensor
products.

\begin{claim}[Pieri formul{\ae}]
\begin{align*}
\Sch\lambda\KK^{\n*}\otimes\Sym[k]\KK^{\n*}&\simeq
\bigoplus_{\substack{\lambda',\ |\lambda'|=|\lambda|+k\\
\lambda'_i\ge\lambda_i\ge\lambda'_{i+1}}}\Sch{\lambda'}\KK^{\n*}
&&\text{as $\GL\n$-modules,}\\
\Sch{\langle\lambda\rangle}\KK^{\n*}\otimes\Sym[k]\KK^{\n*}&\simeq
\bigoplus_{\substack{\lambda',\mu,\ |\mu|=|\lambda|-p\\|\lambda'|=|\mu|+k-p,\ p\le k\\
\lambda_i,\lambda'_i\ge\mu_i\ge\lambda_{i+1},\lambda'_{i+1}}}\Sch{\langle\lambda'\rangle}\KK^{\n*}
&&\text{as $\Sp\n$-modules.}
\end{align*}
\end{claim}

In more visual terms, the Young diagrams of various $\lambda'$ are
obtained from that of $\lambda$ by first removing $p$ boxes from
the right of some rows ($p=0$ for~$\GL\n$) and then adding $k-p$
boxes on the right of some rows in such a way that the horizontal
positions of removed or added boxes do not overlap with those in
other rows and with lower rows of~$\lambda$.

The algebraic study of natural differential operations on tensor
bundles involves polynomial maps between various tensor spaces
which are equivariant with respect to classical groups. These maps
can be described with the aid of classical invariant theory.

\begin{theorem}\label{eq.maps}
Let $G$ be one of the classical linear groups $\GL\n$, $\SL\n$,
$\Or\n$, $\SO\n$, $\Sp\n$ and $V_1,\dots,V_s,W$ be tensor spaces
over~$\KK^{\n}$. Every $G$-equivariant polynomial map
$V_1\times\dots\times V_s\to W$ is obtained by composition and
linear combination from the following basic tensor operations:
\begin{enumerate}

\item\label{tensor} tensor product of elements of various~$V_i$
(maybe occurring repetitively);

\item\label{inv.tens} tensor product with basic $G$-invariant
tensors, which are: the identity operator (\;=~the Kronecker
delta), the co- and contravariant skew-symmetric $\n$-tensors
$\det$ and~$\dual\det$ (for $G=\SL\n,\SO\n$), the co- and
contravariant metric tensors $\metric$ and~$\dual\metric$ (for
$G=\Or\n,\SO\n,\Sp\n$);

\item\label{contr} (partial) contraction;

\item\label{perm} permutation of indices.

\end{enumerate}
\end{theorem}

\begin{proof}
The problem reduces to a description of polynomial functions on
$V_1\times\dots\times V_s\times W^*$ that are linear with respect
to~$W^*$. The assertion in this case stems from the symbolic
method of classical invariant theory \cite[9.5]{inv.th}: it
suffices to involve basic operations (\ref{tensor}),
(\ref{inv.tens}), (\ref{contr}) with basic tensors
$\det,\dual\det,\metric,\dual\metric$ (depending on~$G$).
Contractions involving only indices of $W^*$ correspond to tensor
products with identity operators. Contracting indices of $W^*$
with those of $V_1,\dots,V_s$, and basic tensors in various orders
corresponds to~(\ref{perm}).
\end{proof}

\section{Finiteness theorems}
\label{finite}

\subsection{Reduction to multilinear case}

We start the study of natural differential operations on tensor
fields by some easy reductions.

Let $\V,\W$ be tensor bundles with typical fibres $V,W$. An
algebraic natural differential operation $\DD:\Sec\V\to\Sec\W$ of
order $\le\kk$ is given by a $\JGL[\kk+1]\n$-equivariant
polynomial map $\dd:\JF{V}\to W$. Since $\JGL[\kk+1]\n$ acts on
the vector spaces $\JF{V},W$ linearly, the homogeneous components
of $\dd$ are equivariant maps, too. Thus it suffices to study
homogeneous operations, i.e., those $\DD$ corresponding to
homogeneous~$\dd$.

Assume that $\DD$ is homogeneous of degree $\deg\DD:=\deg\dd=d$.
The polarization of $\dd$ yields a multilinear equivariant map
\begin{equation*}
\dd:\underbrace{\JF{V}\times\dots\times\JF{V}}_{d~\text{times}}\to{W}.
\end{equation*}
(We denote it by the same letter, because the initial homogeneous
map is the restriction of the multilinear map to the diagonal.) By
decomposing the $\GL\n$-modules $V,W$ into irreducibles and by
multilinearity, we reduce $\dd$ to finitely many equivariant
linear maps of the form
\begin{equation}\label{polylin}
\dd:\JF{V_1}\otimes\dots\otimes\JF{V_d}\to W,
\end{equation}
where $V_i,W$ are now assumed to be irreducible $\GL\n$-modules.
In this case we say that the associated tensor bundles $\V_i,\W$
are \emph{indecomposable}. Thus we have reduced (to a certain
extent) the study of arbitrary algebraic natural differential
operations between tensor bundles to the case of multilinear
natural operations on indecomposable tensor bundles.

Consider the map~(\ref{polylin}). Since
\begin{equation*}
\JF{V_1}\otimes\dots\otimes\JF{V_d}=\bigoplus_{l_1,\dots,l_d\le\kk}
\Sym[l_1]\KK^{\n*}\otimes\dots\otimes\Sym[l_d]\KK^{\n*}\otimes
V_1\otimes\dots\otimes V_d,
\end{equation*}
$\dd$~decomposes into a sum of $\GL\n$-equivariant maps
\begin{equation*}
\dd_{l_1,\dots,l_d}:
\Sym[l_1]\KK^{\n*}\otimes\dots\otimes\Sym[l_d]\KK^{\n*}\otimes
V_1\otimes\dots\otimes V_d\to W.
\end{equation*}

For any irreducible $\GL\n$-module
$U=\Sch\lambda\KK^{\n*}\otimes\det^p$ corresponding to a virtual
partition $\overline\lambda=\lambda-(p^{\n})$, put
$|U|=|\overline\lambda|=|\lambda|-\n{p}$.

\begin{lemma}\label{hom.der}
Suppose $\delta_{l_1,\dots,l_d}\ne0$; then
$l_1+\dots+l_d=|W|-\sum|V_i|$. In other words, every natural
multilinear differential operation on indecomposable tensor
bundles is homogeneous with respect to the total order of
derivation.
\end{lemma}

\begin{proof}
Consider the subgroup of homotheties $\KK^{\times}\subset\GL\n$.
For any irreducible $\GL\n$-module $U$ we have: $t\cdot
u=t^{-|U|}u$, $\forall t\in\KK^{\times}$, $u\in U$. Now the lemma
stems from the $\KK^{\times}$-equivariance
of~$\delta_{l_1,\dots,l_d}$.
\end{proof}

\subsection{Finiteness}

Now we prove two general finiteness results. The first one is
easy.

\begin{theorem}
Given two tensor bundles $\V,\W$, the differential order of an
algebraic natural operation $\DD$ from $\V$ to $\W$ is
$O(\deg\DD)$.
\end{theorem}

\begin{corollary}[{cf.~\cite[\S4]{inv.op}}]
Natural algebraic differential operations of degree $\le d$
between two given tensor bundles form a finite-di\-men\-sional
space.
\end{corollary}

\begin{proof}
The polarization and decomposition of $\V,\W$ into indecomposables
reduces the problem to the case where
$\DD:{\Sec{\V_1}\times\dots\times\Sec{\V_d}}\to\Sec\W$ is a
multilinear operation on indecomposable tensor bundles. Now by
Lemma~\ref{hom.der}, the order of $\DD$ is $\max\{l_i\}$ (over
$(l_1,\dots,l_d)$ such that $\dd_{l_1,\dots,l_d}\ne0$)
$\le|W|-\sum|V_i|\le|W|+ d\cdot\max\{-|V_i|\}=O(d)$.
\end{proof}

The second result is much stronger.

\begin{theorem}
Given a tensor bundle $\V$ and two positive integers $\kk,d$,
there exist finitely many indecomposable tensor bundles $\W_i$ and
natural operations $\DD_i:\JF\V\to\W_i$ of degree $\le d$ in
partial derivatives such that every natural operation
$\DD:\JF\V\to\W$ of degree $\le d$ in partial derivatives is
represented as
\begin{equation*}
\DD\vv=\sum_{i,p}\Phi_{ip}(\DD_i\vv\otimes\vv^{\otimes p}),\qquad
\forall\vv\in\Sec\V,
\end{equation*}
where $\Phi_{ip}:\W_i\otimes\Sym[p]\V\to\W$ are natural linear
maps (i.e., compositions of contractions, permutations of indices,
tensor product with the identity operator, and linear
combinations, by Theorem~\ref{eq.maps}).
\end{theorem}

\begin{proof}
Passing to jets, we may reformulate the assertion as follows:
there exist finitely many irreducible $\GL\n$-modules $W_i$ and
$\JGL[\kk+1]\n$-equi\-vari\-ant polynomial maps
\begin{align*}
\dd_i:\JF{V}&=\bigoplus_{l=0}^{\kk}\Sym[l]\KK^{\n*}\otimes V\to W_i\\
\intertext{of degree $\le d$ in the coordinates of}
\JF{V_+}&=\bigoplus_{l=1}^{\kk}\Sym[l]\KK^{\n*}\otimes V
\end{align*}
such that every $\JGL[\kk+1]\n$-equivariant polynomial map
$\dd:\JF{V}\to W$ of degree $\le d$ in $\JF{V_+}$ is of the form
\begin{equation*}
\dd(\vv)=\sum_{i,p}\Phi_{ip}(\dd_i(\vv)\otimes\vv_0^{\otimes p}),
\qquad\forall\vv=\vv_0+\vv_+\in\JF{V}=V\oplus\JF{V_+},
\end{equation*}
where $\Phi_{ip}:W_i\otimes\Sym[p]V\to W$ are
$\GL\n$-equi\-vari\-ant linear maps. The
$\JGL[\kk+1]\n$-equiv\-ari\-ance condition means that the maps are
$\GL\n$-equi\-vari\-ant and $\NGL[\kk+1]\n$-in\-vari\-ant.

All polynomial functions on $\JF{V}$ of degree $\le d$ in
$\JF{V_+}$ form a free $\Sym{V^*}$-module
$\Sym{V^*}\otimes\Mor[d]{\JF{V_+}}\KK$ of finite rank. The
$\NGL[\kk+1]\n$-in\-vari\-ant functions form a $\GL\n$-stable
submodule~$\Inv$. Clearly, operations of order $\le\kk$ and degree
$\le d$ in partial derivatives are identified with $\GL\n$-fixed
elements of $\Inv\otimes W$ or $\GL\n$-equi\-vari\-ant linear maps
$W^*\to\Inv$.

As a submodule of a N{\"o}therian module, $\Inv$~is finitely
generated. Choose finitely many $\dd_i:W_i^*\to\Inv$ whose images
generate $\Inv$ as an $\Sym{V^*}$-module. By complete reducibility
of $\GL\n$-modules, every $\dd:W^*\to\Inv$ lifts to a
$\GL\n$-equi\-vari\-ant linear map
$W^*\to\bigoplus\Sym{V^*}\otimes W_i^*$ along the module
epimorphism $\bigoplus\Sym{V^*}\otimes W_i^*\to\Inv$. The
component mappings $W^*\to\Sym[p]{V^*}\otimes W_i^*$ are nothing
but~$\Phi_{ip}^*$.
\end{proof}

\section{Linear natural operations}
\label{linear}

In this section we deduce the known classification of natural
linear differential operations on tensor bundles in a short and
conceptual way using our algebraic approach (the IT-reduction).
Apparently, the exterior differential is essentially the unique
such operation.

\begin{theorem}\label{lin.op}
Every natural linear differential operation of order $>0$ on
tensor bundles is obtained from the exterior differential by
composition with tensor operations (contraction, permutation of
indices, tensor product with the identity operator) and linear
combination.
\end{theorem}

\begin{remark}
This theorem is usually referred to as the Schouten theorem,
although Schouten just formulated it in 1951 without proof. It was
proved for differential forms by Palais (1959), for arbitrary
covariant tensors by Leicher (1973), and in full generality by
Rudakov (1974), Chuu-Lian Terng (1976), and Kirillov (1977),
see~\cite{inv.op}. For operations of order~1, a proof based on the
IT-reduction was first obtained by Smirnov~\cite{diplom}.
\end{remark}

\begin{proof}
Let $\V,\W$ be two tensor bundles with typical fibres $V,W$. A
linear differential operation $\DD:\JF\V\to\W$ is determined by a
$\GL\n$-equi\-vari\-ant linear map $\dd:\JF{V}\to W$ which
vanishes on $\ngl[\kk+1]\n\cdot\JF{V}$. Without loss of generality
we may assume that $\V,\W$ are indecomposable. (Injections and
projections onto indecomposable summands are given by tensor
operations.)

We prove that $\Sym[\kk]\KK^{\n*}\otimes V\subseteq\Ker\dd$ unless
$\kk=1$, $V\simeq\Ext[m]\KK^{\n*}$, $m<\n$. It suffices to show
that $\dg\cdot\vv$ span $\Sym[\kk]\KK^{\n*}\otimes V$ for some
$\dg\in\ngl[\kk+1]\n$, $\vv\in\JF{V}$. Suppose that
$V=\Sch\lambda\KK^{\n*}\otimes\det^d$ corresponds to a virtual
partition $\overline\lambda=\lambda-(d^{\n})$.

If $\kk>1$, then $\Sym[\kk]\KK^{\n*}\otimes V$ is spanned by
$\dg\cdot\vv$, $\dg=\dg[\kk]\in\Sym[\kk]\KK^{\n*}\otimes\KK^{\n}$,
$\vv=\vv[1]\in\KK^{\n*}\otimes V$. Indeed, take
\begin{equation*}
\dg=(\x{i})^{\kk}\otimes\e{i},\qquad
\vv=\x{i}\otimes\vv[\lambda]^-,
\end{equation*}
where $\e1,\dots,\e\n$ are the standard basic vectors
of~$\KK^{\n}$. (Here and below in the proof, we do \emph{not} sum
over~$i$.) By (\ref{Lie.act}) we have
\begin{equation*}
\dg\cdot\vv=-(\kk\overline\lambda_i+1)(\x{i})^{\kk}\otimes\vv[\lambda]^-.
\end{equation*}
For $i=\n$, $\dg\cdot\vv$~generates $\Sym[\kk]\KK^{\n*}\otimes V$
as a $\GL\n$-module, because it is the product of a highest and a
lowest weight vectors.

If $\kk=1$, then we are left with
$\dg\in\Sym[2]\KK^{\n*}\otimes\KK^{\n}$, $\vv\in V$. Take
\begin{equation*}
\dg=(\x{i})^2\otimes\e{i},\ \vv=\vv[\lambda]^- \implies
\dg\cdot\vv=-2\overline\lambda_i\x{i}\otimes\vv[\lambda]^-.
\end{equation*}
If $\overline\lambda_{\n}\ne0$, then we conclude as above.
Otherwise suppose $\overline\lambda_1>1$. Put
\begin{gather*}
\dg'=\x{i}\x{j}\otimes\e{i}\implies
\begin{cases}
\dg'\cdot\vv=-\overline\lambda_i\x{j}\otimes\vv-\x{i}\otimes\vv',&\\
\dg\cdot\vv'=-2(\overline\lambda_i-1)\x{i}\otimes\vv',&\text{where}
\end{cases}\\
\vv'=\sum_{p=1}^{\lambda_i}
\Alt[\lambda]\left((\x1)^{\otimes\lambda_1}\cdots
(\x{i})^{\otimes{p-1}}\otimes\x{j}\otimes(\x{i})^{\otimes\lambda_i-p}
\cdots(\x\n)^{\otimes\lambda_{\n}}\right).
\end{gather*}
For $i=1$, $j=\n$, we obtain $\x\n\otimes\vv[\lambda]^-\in\Ker\dd$
and conclude as above.

We are left with the case $\kk=1$, $V=\Ext[m]\KK^{\n*}$, $m<\n$.
Here
\begin{equation*}
\KK^{\n*}\otimes V=
\Ext[m+1]\KK^{\n*}\oplus\Sch{(2,1^{m-1})}\KK^{\n*}.
\end{equation*}
All $\dg\cdot\vv$ are killed by $\Alt[(1^{m+1})]$, hence span the
2-nd direct summand. The unique natural operation, given by the
projection $\Alt[(1^{m+1})]$ onto the 1-st summand, is the
exterior differential, cf.~Example~\ref{ext.diff}.
\end{proof}

\begin{remark}
It would be interesting to reproduce by this method the
classification of natural bilinear differential operators obtained
by Grozman~\cite{bilin}. However the computations here will be
more involved.
\end{remark}

\section{Natural operations on symplectic manifolds}
\label{symplectic}

\subsection{Geometric structures}
\label{geom.str}

Manifolds are often equipped with additional geometric structures,
and it is important to study differential operations that are
``natural'' with respect to these structures. There are several
possible ways to formalize these concepts. Here we adopt the
following one.

\begin{definition}
A \emph{type of geometric structures} is a functor on $\Man\n$
associating with each $\n$-manifold $\M$ a subsheaf $\GS=\GS(\M)$
in the sheaf of sections of a natural bundle $\F\to\M$ so that,
for any open embedding $\M'\hookrightarrow\M$, $\GS(\M')$~is the
pullback of~$\GS(\M)$. A \emph{geometric structure} (of given
type) is a section $\gs\in\Sec\GS$.
\end{definition}

\begin{remark}
As a rule, $\GS$~consists of sections satisfying a certain natural
differential equation.
\end{remark}

\begin{example}
A Riemannian structure is given by a section of a natural bundle
$(\Sym[2]\T*)^+$ of positive quadratic forms on tangent spaces
($\KK=\RR$).
\end{example}

\begin{example}
A symplectic structure is a section of
\begin{equation*}
\GS=\{\omega\in\Sec{(\DF[2])^{\reg}}\mid d\omega=0\},
\end{equation*}
where $(\DF[2])^{\reg}$ is a natural bundle of non-degenerate
2-forms.
\end{example}

\begin{example}\label{Poisson}
A Poisson structure is a section of
\begin{equation*}
\GS=\{\beta\in\Sec{\Ext[2]\T{}}\mid[\beta,\beta]=0\},
\end{equation*}
where
$[\beta,\beta]=\beta^{ij}\partial_i\beta^{kl}\;\partial_j\wedge\partial_k\wedge\partial_l$
is the Schouten bracket.
\end{example}

\begin{definition}\label{nat.op&str}
A \emph{natural differential operation} on manifolds with a
geometric structure of type $\GS$ is a natural operation (in the
sense of~\ref{nat.op})
\begin{equation*}
\DD:\Sec\V\times\Sec\GS\to\Sec\W,
\end{equation*}
where $\V,\W$ are natural bundles.
\end{definition}

In other words, $\DD\vv$~depends on $\vv$ and $\gs$ (as a
parameter), and on their partial derivatives (up to a certain
order) in a way independent of chosen local coordinates, where
$\vv$ and $\gs$ are (local) sections of $\V$ and~$\GS$,
respectively.

\begin{example}
The covariant differential on tensor fields is a natural operation
on Riemannian manifolds. The Poisson bracket of functions is a
natural operation on symplectic or Poisson manifolds.
\end{example}

Alike ordinary natural differential operations, those on manifolds
with an additional structure are uniquely determined by their
action on jets:
\begin{equation*}
\dd:\JF{V}\times\JF[l]\Gs\to W,
\end{equation*}
where $\JF[l]\Gs\subseteq\JF[l]F$ is the space of jets (at some
point) of local sections from~$\GS$, and $\dd$ is
$\JGL[r]\n$-equivariant for an appropriate~$r$.

In the presence of an additional geometric structure, the group of
coordinate transformations, which must preserve the coordinate
expression of a natural operation, often can be reduced by the
following standard invariant-theoretic trick.

\begin{definition}
Suppose a Lie (or algebraic) group $G$ acts on a manifold
(algebraic variety) $X$ and $H\subset G$ is a Lie (algebraic)
subgroup. An $H$-stable submanifold (subvariety) $Y\subset X$ is
said to be a \emph{$(G,H)$-section} if the natural map
$G\gtimes{H}Y\to X$ is an isomorphism.
\end{definition}

The significance of this notion is that $G$-invariant maps of $X$
are in a bijective correspondence (by restriction to~$Y$) with
$H$-invariant maps of~$Y$.

In particular, for some important geometric structures it happens
that $\JF[l]\Gs$ admits a
$\bigl(\JGL[r]\n,\JF[r]{G}\bigr)$-section $\JF[l]{\Gs_0}$ for some
subgroup $\JF[r]{G}\subset\JGL[r]\n$. Therefore natural
differential operations are in a bijective correspondence with
$\JF[r]{G}$-invariant maps
\begin{equation*}
\dd:\JF{V}\times\JF[l]{\Gs_0}\to W.
\end{equation*}

\begin{example}[{\cite[\S3]{curv}}]
The space of $\kk$-jets of Riemannian metrics
\begin{align*}
\JF\Gs&=(\Sym[2]\RR^{\n*})^{+}\oplus
\bigoplus_{l=1}^{\kk}\Sym[l]\RR^{\n*}\otimes\Sym[2]\RR^{\n*}\\
\intertext{admits a $\bigl(\JGL[\kk+1]\n,\GL\n\bigr)$-section}
\JF{\Gs_0}&=(\Sym[2]\RR^{\n*})^{+}\oplus
\bigoplus_{l=2}^{\kk}\Ker\Symm[l+1],
\end{align*}
where
$\Symm[l+1]:\Sym[l]\RR^{\n*}\otimes\Sym[2]\RR^{\n*}\to\Sym[l+1]\RR^{\n*}\otimes\RR^{\n*}$
is the symmetrization in the first $l+1$ indices. In fact,
$\Ker\Symm[l+1]\simeq\Sch{(l,2)}\RR^{\n*}$. Moving a metric to
$\JF{\Gs_0}$ by a coordinate transformation corresponds to writing
the metric in geodesic coordinates with center at given point.
\end{example}


\subsection{Symplectic structure}

Jets of symplectic structures form a homogeneous space, i.e., one
may take just one point for a section, by the Darboux theorem. It
is instructive to give a purely algebraic proof of this fact in
the spirit of this paper.

\begin{claim}[Formal Poincar\'e Lemma]
The De~Rham complex
\begin{equation*}
\cdots\overset{d}\longrightarrow
\Sym[l]\KK^{\n*}\otimes\Ext[m]\KK^{\n*}\overset{d}\longrightarrow
\Sym[l-1]\KK^{\n*}\otimes\Ext[m+1]\KK^{\n*}\overset{d}\longrightarrow\cdots
\end{equation*}
is exact.
\end{claim}

\begin{proof}
By the Pieri formul{\ae}, we have
\begin{equation*}
\Sym[l]\KK^{\n*}\otimes\Ext[m]\KK^{\n*}\simeq
\Sch{(l+1,1^{m-1})}\KK^{\n*}\oplus\Sch{(l,1^m)}\KK^{\n*}.
\end{equation*}
(In fact,
$\Sch{(l+1,1^{m-1})}\KK^{\n*}\subset\Sym[l]\KK^{\n*}\otimes\Ext[m]\KK^{\n*}$
and $\Sch{(l,1^m)}\KK^{\n*}$ embeds as $\Ker\Symm[l+1]$.) The map
$d:\Sym[l]\KK^{\n*}\otimes\Ext[m]\KK^{\n*}\to\Sym[l-1]\KK^{\n*}\otimes\Ext[m+1]\KK^{\n*}$
is nothing but the alternation in the last $m+1$ indices. It maps
$\Sch{(l,1^m)}\KK^{\n*}$ isomorphically into
$\Sym[l-1]\KK^{\n*}\otimes\Ext[m+1]\KK^{\n*}$, while its kernel
$\Sch{(l+1,1^{m-1})}\KK^{\n*}$ is exactly the image of
$\Sym[l+1]\KK^{\n*}\otimes\Ext[m-1]\KK^{\n*}$.
\end{proof}

\begin{claim}[Formal Darboux Theorem]
Every jet
$\omega(\x{})=\omega_0+\omega_1(\x{})+\dots+\omega_k(\x{},\dots,\x{})\in\JP\n\otimes\Ext[2]\KK^{\n*}$
such that $\omega_0$ is non-degenerate and $d\omega=0$ is
$\JGL[\kk+1]\n$-equivalent to~$\omega_0$.
\end{claim}

\begin{proof}
We shall successively kill the non-constant terms of $\omega$ by
transformations $\x{}\mapsto\g(\x{})=
\x{}+\g[2](\x{},\x{})+\dots+\g[\kk+1](\x{},\dots,\x{})$
from~$\NGL[\kk+1]\n$. Arguing by induction on~$\kk$, we may assume
that $\omega_l=0$, $0<l<\kk$. Take $\g\in\NGL[\kk+1]\n$ such that
$\g[l]=0$, $2\le l\le\kk$. Using Lemma~\ref{tensor.jet} and its
corollary and adopting the notation therein, we have
\begin{multline*}
(\g\cdot\omega)(\x{})
=\omega(\x{})+(\kk+1)\dR\bigl(\g[\kk+1](\x{},\dots,\x{},\cdot)\bigr)\omega_0\\
=\omega(\x{})+(\kk+1)\act(\g[\kk+1]\otimes\omega_0)
=\omega(\x{})+d(\g[\kk+1]*\omega_0),
\end{multline*}
where $(\cdot)*\omega_0$ is the lowering of the upper index by
contraction with $\omega_0$ in its 2-nd index. Thus $\omega_{\kk}$
may be shifted by an arbitrary vector in the image of
\begin{equation*}
\Sym[\kk+1]\KK^{\n*}\otimes\KK^{\n}
\overset{*\omega_0}{\relbar\joinrel\longrightarrow}
\Sym[\kk+1]\KK^{\n*}\otimes\KK^{\n*}\overset{d}\longrightarrow
\Sym[\kk]\KK^{\n*}\otimes\Ext[2]\KK^{\n*}.
\end{equation*}
By the formal Poincar\'e lemma, the image consists of all (jets
of) closed forms. Hence $\omega_{\kk}$ may be shifted to~0.
\end{proof}

Acting by $\GL\n$ we may obtain that $\omega_0$ be the standard
symplectic form on $\KK^{\n}$ such that
\begin{equation*}
\omega_0(\e{i},\e{j})=
\begin{cases}
 1,& i+j=\n+1,\ i<j, \\
-1,& i+j=\n+1,\ i>j, \\
 0,& \text{otherwise.}
\end{cases}
\end{equation*}
We identify $\KK^{\n}$ with $\KK^{\n*}$ by lowering the indices.
The symplectic form $\omega_0^*$ on $\KK^{\n*}$ induced from
$\omega_0$ is given by
$\omega_0^*(\x{i},\x{j})=\omega_0(\e{i},\e{j})$.

The stabilizer of $\omega_0$ in $\JGL\n$ is the group $\JSp\n$ of
$\kk$-jets of symplectomorphisms $\KK^{\n}\to\KK^{\n}$ at~$0$. It
has Levi decomposition $\JSp\n=\Sp\n\rightthreetimes\NSp\n$, with
the unipotent radical $\NSp\n=\NGL\n\cap\JSp\n$. Let us describe
the Lie algebra $\jsp\n$ of~$\JSp\n$.

\begin{lemma}\label{jsp}
$\dg\in\jsp\n\iff\dg[l]*\omega_0\in\Sym[l+1]\KK^{\n*},\ \forall
l=1,\dots,\kk$ (in particular, $\dg[1]\in\spl\n$)
\end{lemma}

\begin{proof}
By the corollary of Lemma~\ref{tensor.jet},
\begin{equation*}
\dg\cdot\omega_0=d(\dg*\omega_0)=\sum_{l=1}^{\kk}d(\dg[l]*\omega_0)=0\iff
d(\dg[l]*\omega_0)=0,\ \forall l=1,\dots,\kk.
\end{equation*}
We have $\dg[l]*\omega_0\in\Sym[l]\KK^{\n*}\otimes\KK^{\n*}\simeq
\Sym[l+1]\KK^{\n*}\oplus\Sch{(l,1)}\KK^{\n*}$, and the 1-st
summand is exactly~$\Ker{d}$.
\end{proof}

By the above reasoning, natural differential operations on
symplectic manifolds are given by $\JSp[r]\n$-equivariant maps
$\dd:\JF{V}\to W$.

\begin{remark}
Another approach to natural differential operations on symplectic
manifolds is to extend the notion of geometric objects by
considering fibre bundles $\F=\SRep\M\gtimes{\JSp\n}F$ associated
with the symplectic coframe bundle $\SRep\M\to\M$ of order~$\kk$.
The latter consists of $\kk$-jets of symplectic coordinate systems
transforming the symplectic form $\omega$ on $\M$ into the
standard symplectic form $\omega_0$ on~$\KK^{\n}$. Note that $F$
is acted on only by $\JSp\n$, not~$\JGL\n$, but if the action
extends to~$\JGL\n$, then
$\F\simeq\SRep\M\gtimes{\JSp\n}\JGL\n\gtimes{\JGL\n}F\simeq\Rep\M\gtimes{\JGL\n}F$
is a natural bundle in the sense of~\ref{geom.quant}.

Now natural differential operations $\DD:\Sec\V\to\Sec\W$ are
defined as those having one and the same expression in all
\emph{symplectic} coordinate systems. This is equivalent to the
map of jets $\dd:\JF{V}\to W$ be $\JSp[r]\n$-equi\-vari\-ant.
\end{remark}

\subsection{Linear operations}

Now we describe linear natural differential operations on
symplectic manifolds in a way similar to Section~\ref{linear}. The
classification of such operations, straightforward via our
approach, was first obtained by Rudakov \cite{irrep(ham)} in a
rather indirect way from the study of irreducible representations
of certain infinite-dimensional Lie algebras.

On symplectic $\n$-manifolds, the natural isomorphism
$\DF[m]\simeq\DF[m*]$ allows to define the differential operation
$d^*:\Sec{\DF[m]}\to\Sec{\DF[m-1]}$ contragredient to
$d:\Sec{\DF[m-1]}\to\Sec{\DF[m]}$. Actually $d^*$ is the
composition of $d:\Sec{\DF[\n-m]}\to\Sec{\DF[\n-m+1]}$ with the
natural isomorphisms
$\DF[m]\simeq\DF[m*]\stackrel\sim\longrightarrow\DF[\n-m]$ and
$\DF[\n-m+1]\stackrel\sim\longrightarrow(\DF[m-1])^*\simeq\DF[m-1]$
given by tensor operations.

\begin{theorem}
Every natural linear differential operation of order $>0$ on
tensor bundles over symplectic manifolds is obtained from the
exterior differential $d$ and the symplectic Laplacian $dd^*$ by
composition with tensor operations (contraction, permutation of
indices, tensor product with the symplectic form $\omega$ or the
dual bivector~$\omega^*$) and linear combination.
\end{theorem}

\begin{proof}
Similarly to the proof of Theorem~\ref{lin.op}, it suffices to
classify $\Sp\n$-equi\-vari\-ant linear maps $\dd:\JF{V}\to W$
which vanish on $\nsp[\kk+1]\n\cdot\JF{V}$, where $V,W$ are
irreducible $\Sp\n$-modules. Assuming
$V=\Sch{\langle\lambda\rangle}\KK^{\n*}$, we prove that
$\Sym[\kk]\KK^{\n*}\otimes V\subseteq\nsp[\kk+1]\n\cdot\JF{V}
\subseteq\Ker\dd$ unless $\lambda=(1^m)$, $\kk\le2$.

Recall that, by Lemma~\ref{jsp}, we may identify
$\bigoplus_{l=2}^{\kk+1}\Sym[l+1]\KK^{\n*}$ with $\nsp[\kk+1]\n$
by raising an index by contraction with $\omega_0^*$ in its 1-st
index.

If $\kk>1,\lambda_1\ne1$, then $\Sym[\kk]\KK^{\n*}\otimes V$ is
spanned by $\Sym[\kk+1]\KK^{\n*}\cdot(\KK^{\n*}\otimes V)$.
Indeed, take
\begin{align*}
\dg&=(\x\n)^{\kk+1},&\vv&=\x1\otimes\vv[\lambda]^-,\\
\dg'&=(\kk+1)\x1(\x\n)^{\kk},&\vv'&=\x\n\otimes\vv[\lambda]^-;
\end{align*}
then by~(\ref{Lie.act}) we have
\begin{gather*}
\begin{aligned}
\dg\cdot\vv&=(\x\n)^{\kk}\otimes\vv[\lambda]^-+\kk\x1(\x\n)^{\kk-1}\otimes\bar\vv,\\
\dg'\cdot\vv'&=(\kk\lambda_1-1)(\x\n)^{\kk}\otimes\vv[\lambda]^-+\kk(\kk-1)\x1(\x\n)^{\kk-1}\otimes\bar\vv,
\qquad\text{where}
\end{aligned}\\
\bar\vv=\sum_{p=1}^{\lambda_1}
\Alt[\lambda]\left((\x1)^{\otimes{p-1}}\otimes\x\n\otimes(\x1)^{\otimes\lambda_1-p}
\otimes(\x2)^{\otimes\lambda_2}\cdots
\otimes(\x{\n/2})^{\otimes\lambda_{\n/2}}\right).
\end{gather*}
We obtain $\kk(\lambda_1-1)\x\n\otimes\vv[\lambda]^-\in
\Sym[\kk+1]\KK^{\n*}\cdot(\KK^{\n*}\otimes V)$. Since
$\lambda_1\ne1$, this vector generates $\Sym[\kk]\KK^{\n*}\otimes
V$ as an $\Sp\n$-module, because it is the product of a highest
and a lowest weight vectors.

In the case $\kk=1$ or $\lambda_1=1$, take
\begin{align*}
\dg&=3\x1(\x\n)^2,&\vv&=(\x\n)^{\kk-1}\otimes\vv[\lambda]^-,\\
\dg'&=3(\x1)^2\x\n,&\vv'&=(\x\n)^{\kk-1}\otimes\bar\vv;
\end{align*}
then
\begin{align*}
\dg\cdot\vv&=(2\lambda_1+1-\kk)(\x\n)^{\kk}\otimes\vv[\lambda]^-+2\x1(\x\n)^{\kk-1}\otimes\bar\vv,\\
\dg'\cdot\vv'&=-2\lambda_1(\x\n)^{\kk}\otimes\vv[\lambda]^-+2(\lambda_1-1-\kk)\x1(\x\n)^{\kk-1}\otimes\bar\vv.
\end{align*}
We derive that
$(\lambda_1+1-\kk)(2\lambda_1-1-\kk)(\x\n)^{\kk}\otimes\vv[\lambda]^-\in
\Sym[3]\KK^{\n*}\cdot{(\Sym[\kk-1]\KK^{\n*}\otimes V)}$. As above,
this vector is nonzero and generates $\Sym[\kk]\KK^{\n*}\otimes V$
as an $\Sp\n$-module unless $\kk\le2$, $\lambda_1\le1$.

We are left with the case $\lambda=(1^m)$, $\kk\le2$. For $\kk=1$,
we have
\begin{equation*}
\KK^{\n*}\otimes V\simeq\Sch{\langle1^{m-1}\rangle}\KK^{\n*}\oplus
\Sch{\langle1^{m+1}\rangle}\KK^{\n*}\oplus
\Sch{\langle2,1^{m-1}\rangle}\KK^{\n*}.
\end{equation*}
The projections onto the first two summands correspond to $d^*$
and~$d$, respectively, while the 3-rd summand is exactly
$\Sym[3]\KK^{\n*}\cdot V$.

For $\kk=2$, we have
\begin{equation*}
\Sym[2]\KK^{\n*}\otimes V\simeq
\Sch{\langle1^m\rangle}\KK^{\n*}\oplus
\Sch{\langle2,1^{m-2}\rangle}\KK^{\n*}\oplus
\Sch{\langle2,1^m\rangle}\KK^{\n*}\oplus
\Sch{\langle3,1^{m-1}\rangle}\KK^{\n*}.
\end{equation*}
It is easy to verify that the last 3 summands are in
$\Sym[3]\KK^{\n*}\cdot(\KK^{\n*}\otimes V)$. Indeed, for $i<\n$ we
have
\begin{equation*}
(\x1)^3\cdot\bigl(\x{i}\otimes(\x\n\wedge\x2\wedge\dots\wedge\x{m})\bigr)=
-2\x1\x{i}\otimes(\x1\wedge\dots\wedge\x{m})=:w.
\end{equation*}
For $i=1$, $w$~is a lowest weight vector
of~$\Sch{\langle3,1^{m-1}\rangle}\KK^{\n*}$. For ${i=m+1}$, the
invariant projector $\Alt[(2,1^m)]$ maps $w$ to a lowest weight
vector
${(-1)^{m+1}\x1\otimes\x1\wedge\dots\wedge\x{m+1}}\in\Sch{\langle2,1^m\rangle}\KK^{\n*}$.
For $i=\n+1-m$, another invariant projector, namely the
contraction with $\omega_0^*$ in, say, the first and last indices,
maps $w$ to a lowest weight vector
${\x1\otimes\x1\wedge\dots\wedge\x{m-1}}\in\Sch{\langle2,1^{m-2}\rangle}\KK^{\n*}$.
Since these irreducible modules occur in $\Sym[2]\KK^{\n*}\otimes
V$ exactly once, they are in
$\Sym[3]\KK^{\n*}\cdot(\KK^{\n*}\otimes V)$ by the Schur lemma.
The projection onto the remaining 1-st summand corresponds
to~$dd^*$.
\end{proof}

\begin{remark}
Natural bilinear differential operators on symplectic manifolds
were partially classified by Grozman~\cite{Sp-bilin}. It would be
interesting to complete the classification using our methods.
\end{remark}

\section{Deformation quantization}
\label{def.quant}

In this section, we address the problem of existence of a natural
deformation quantization on Poisson manifolds. Let us recall the
notion of deformation quantization.

Given a Poisson manifold $\M$ with a Poisson bivector
$\beta\in\Sec[\M]{\Ext[2]\T{}}$ (cf.~Example~\ref{Poisson}), the
sheaf $\OO$ of differentiable functions on $\M$ comes equipped
with the Poisson bracket
$\{f,g\}=\beta(df,dg)=\beta^{ij}\,\partial_if\,\partial_jg$. A
\emph{deformation quantization} is an associative product $\star$
on the sheaf $\OO\[\h\]$ of formal power series with coefficients
in $\OO$ that is $\KK\[\h\]$-linear with respect to infinite
formal sums and is defined on $\OO$ by a formula
\begin{equation}\label{star}
f\star g=fg+\h\{f,g\}+\dots+\h^m\beta_m(f,g)+\cdots,
\end{equation}
where $\beta_m$ ($m=1,2,\dots$) are bilinear differential
operators. The $\star$-prod\-uct may be considered as a
non-commutative deformation of the usual commutative product of
functions, $\h$~being the parameter of deformation, so that the
1-st order term of $f\star g-g\star f$ is $2\{f,g\}$. Deformation
quantization is one of possible approaches to mathematical
foundations of quantum mechanics~\cite{def&quant}.

Several constructions of a deformation quantization for a given
Poisson structure are known: by Moyal, by
De~Wilde--Lecomte~\cite{*prod}, by Fedosov \cite{F-quant},
\cite{def.quant}, by Kontsevich~\cite{quant.Poisson}, etc. But all
of them involve some additional geometric structure on a Poisson
manifold: affine structure, linear connection, etc. A natural
question arises: does there exist a canonical deformation
quantization given by one and the same universal formula for all
Poisson manifolds? In terms of this paper, this may be
reformulated as follows: does there exist a $\star$-product
(\ref{star}) such that its terms $\beta_m$ are natural
differential operations on Poisson manifolds in the sense of
Definition~\ref{nat.op&str}? We answer this question negatively.

\begin{theorem}\label{nat.quant}
There exists no natural deformation quantization on Poisson
manifolds.
\end{theorem}

\begin{proof}
First we prove the theorem for symplectic manifolds. In the
symplectic case, $\beta=\omega^*$ is the bivector dual to the
symplectic form~$\omega$. The operations $\beta_m$ are determined
by the respective bilinear maps of jets
\begin{equation*}
\beta_m:\JP\n\otimes\JP[l]\n\to\KK,
\end{equation*}
which must be  $\JSp[r]\n$-equivariant, $r\ge\max\{\kk,l\}$.

Since $\Sym[\kk]\KK^{\n*}$ are pairwise distinct self-dual
irreducible $\Sp\n$-modules, a non-zero $\Sp\n$-invariant linear
map $\Sym[\kk]\KK^{\n*}\otimes\Sym[l]\KK^{\n*}\to\KK$ exists (and
is then given by the full contraction with
$(\omega_0^*)^{\otimes\kk}$, up to proportionality) iff $\kk=l$.

However, if $\kk>1$, then the unique trivial irreducible factor
$\KK\hookrightarrow\Sym[\kk]\KK^{\n*}\otimes\Sym[\kk]\KK^{\n*}$ is
in~$\Ker\beta_m$. Indeed, consider the action
$\nsp[r]\n:{\JP\n\otimes\JP[l]\n}$. In particular, we have
\begin{equation*}
\Sym[3]\KK^{\n*}\cdot(\Sym[\kk-1]\KK^{\n*}\otimes\Sym[\kk]\KK^{\n*})\subseteq
(\Sym[\kk]\KK^{\n*}\otimes\Sym[\kk]\KK^{\n*})\oplus
(\Sym[\kk-1]\KK^{\n*}\otimes\Sym[\kk+1]\KK^{\n*}).
\end{equation*}
Put $\dg=3(\x1)^2\x{\n}$,
$\vv=(\x1)^{\kk-1}\otimes(\x{\n})^{\kk}$; then
\begin{equation*}
\dg\cdot\vv=(\kk-1)(\x1)^{\kk}\otimes(\x{\n})^{\kk}-2\kk(\x1)^{\kk-1}\otimes\x1(\x{\n})^{\kk}
\end{equation*}
has a non-zero projection to~$\KK$, whence
$\KK\subset\Sym[3]\KK^{\n*}\cdot(\Sym[\kk-1]\KK^{\n*}\otimes\Sym[\kk]\KK^{\n*})$.
Since $\beta_m$ is $\NSp[r]\n$-invariant,
$\Ker\beta_m\supseteq\nsp[r]\n\cdot\bigl(\JP\n\otimes\JP[l]\n\bigr)$,
whence the claim.

It follows that $\beta_m$ is a linear combination of the usual
multiplication $\KK\otimes\KK\to\KK$ and the Poisson bracket
$\KK^{\n*}\otimes\KK^{\n*}\to\KK$. However, it is easy to see that
these two operations cannot be combined into an associative
$\star$-product: in fact, the associativity will be violated
already at order~2 (i.e., in the coefficient at~$\h^2$).

Now it is obvious that a natural quantization cannot exist on
Poisson manifolds of even dimension, because otherwise it would
restrict to a natural quantization of symplectic manifolds as a
particular case.

On a Poisson manifold $\M$ of odd dimension, the Poisson structure
$\beta$ is always degenerate. However, an open subset of $\M$
where $\rk\beta=\max$ has a foliation with symplectic leaves. If
we trivialize this foliation in a neighborhood of a point and
consider functions depending only on the coordinates along the
leaves and coordinate transformations which do not involve the
remaining coordinates, then a natural quantization on $\M$ would
restrict to a natural quantization of the symplectic leaves, a
contradiction.
\end{proof}

\begin{remark}
Theorem~\ref{nat.quant} is a deformation quantization analogue of
the Van~Hove theorem in the theory of geometric quantization
\cite[5.2.2]{geom.quant}, cf.~\cite{no.inv.quant}.
\end{remark}


\end{document}